\documentclass[12pt]{article}
\usepackage{amsmath,amsfonts,amsthm,bm}

\newtheorem{theorem}{\bf Theorem}[section]
\newtheorem{lemma}{\bf Lemma}[section]

\newcommand{\eq}[1]{\label{eq:{#1}}}
\newcommand{\re}[1]{(\ref{eq:{#1}})}
\newcommand{\ct}[1]{\cite{kn:{#1}}}
\newcommand{\by}[1]{\bibitem{kn:{#1}}}
\newcommand{\RR}{\mathbb R}

\newcommand{\dbx}{ \, d{\bf x}\,}
\newcommand{\dbxdt}{ \, d{\bf x}\, dt}
\newcommand\be{\begin{equation}}
\newcommand\ee{\end{equation}}
\newcommand\ba{\begin{eqnarray}}
\newcommand\ea{\end{eqnarray}}
\newcommand{\nn}{\nonumber}

\newcommand{\uu}{{\bf u}}
\newcommand{\vv}{{\bf v}}
\newcommand{\xx}{{\bf x}}
\newcommand{\DD}{\Delta}

\newcommand{\df}{\equiv}

\newcommand\eps{\varepsilon}
\newcommand\ga{\gamma}
\newcommand\gi{\gamma^{-1}}
\newcommand\gd{\dot{\gamma}}
\newcommand\gdd{\ddot{\gamma}}

\newcommand\dt{\partial_{t}}
\newcommand\dx{\partial_{x}}
\newcommand\dy{\partial_{y}}
\newcommand\dxx{\partial_{x}^{2}}
\newcommand\dyy{\partial_{y}^{2}}
\newcommand\dxxx{\partial_{x}^{3}}
\newcommand\dxi{\partial_{x}^{-1}}
\newcommand{\dxii}{\partial_x^{-2}}
\newcommand\dxg{(\partial_{x})_{\ga}}
\newcommand\dxig{(\partial_{x}^{-1})_{\ga}}

\newcommand\dxiidxg{(\dxii\dx)_{\ga}}

\newcommand\sob{H^{s}}

\newcommand\tor{\mathbb{T}}

\newcommand\reel{\mathbb{R}}

\newcommand\D{\mathcal{D}}

\newcommand\ctt{\cos \theta (t)}
\newcommand\stt{\sin \theta (t)}

\title{Global Existence of Infinite Energy Solutions for a Perfect Incompressible Fluid}
\author{{\em Ralph Saxton\thanks{R. Saxton would like to thank Barbara Keyfitz, colleagues and  staff for their support at the Fields Instute, University of Toronto, where  part of this work was carried out.} and Feride T\i\u{g}lay\thanks{F. T\i\u{g}lay would like to thank the Department of Mathematics at the University of Notre Dame for its hospitality.}}\\
{\em Department of Mathematics}\\
{\em University of New Orleans}\\
{\em New Orleans, LA 70148, USA}}
	
\begin{document}
\maketitle

\section{Abstract}
\label{abstract}

This paper provides results on local  and global existence for a class of solutions to the  Euler equations for an incompressible, inviscid fluid. By considering a class of solutions which exhibits a characteristic growth at infinity we obtain an initial value problem for a nonlocal equation.
We establish local well-posedness  in all dimensions and persistence in time of these solutions for three and higher dimensions. We also examine a weaker class of global solutions. 

\section{Introduction}
\label{introduction}
\setcounter{equation}{0}

A fundamental question in the study of fluids concerns the possibility of finite time blow up of solutions to the Euler equations for a perfect, incompressible fluid. It is well known that blow up cannot take place in the two-dimensional case, for solutions defined over a bounded domain  subject to Dirichlet boundary conditions (see for example the work of Wolibner \ct{Wo} and Ebin \ct{Eb}), since smooth data in this case lead to solutions remaining smooth for all time. However the  question  remains open in higher dimensions.

A separate class of solutions  consists of those having   \lq stagnation-point\rq\, form, which  attracted early attention by Weyl \ct{W} and Lin \ct{L},  and provides a set of equations which depend on only a single spatial variable and time.   The  resulting equations, once solved, provide  exact solutions to the full Euler equations. However, the associated growth of the full solutions in certain directions means that the  flows possess, at best, only locally finite kinetic energy. Nevertheless, one may discuss such questions  as finite time blow up for this class, and it  has been shown by Stuart in \ct{S} that this can  take place when  the reduced equations are defined over the real line and those solutions decay at infinity. In this case the corresponding spatial domain for the full equations is $\reel^n$ for $n=2$ and $n=3$, with the  full set of solutions  growing linearly in  the other direction(s).

The evolution of two-dimensional solutions, which can blow up on the unbounded domain $\reel^2$, therefore differs significantly from the  class of globally defined solutions  which exists for bounded subdomains of $\reel^2$. So the consequences of higher-dimensional stagnation-point solutions blowing up might be thought to  result simply from their behavior at infinity rather than bearing on the  question of singularity formation. This view is in some sense strengthened by the results of Childress et al in \ct{CISY} (see also a related result by Cox in \ct{Cox}) where  solutions defined over a two-dimensional, infinite strip were examined. Since blow up was still found over this, smaller, semi-infinite domain, stagnation-point solutions defined over such domains would give the appearance of behaving, generally, much as those defined on the full space. 

A similar approach is implemented by Constantin in \ct{Con} to reduce the Euler equations, periodic in two directions, to a nonlocal Riccati equation and prove the blow-up in finite time by solving these equations on characteristics. 

In this paper we consider a stagnation-point class of solutions defined over $\reel^{n}$ which is spatially periodic in one coordinate direction. In two dimensions, the equations reduce to those of \ct{CISY} and, although we examine slightly different boundary conditions, the same blow up results essentially apply. In three and higher dimensions we however find  a \lq regularizing\rq\, effect not present in solutions which  decay to zero at infinity in the same coordinate direction, and this leads to the existence for all time of all such solutions, stemming from sufficiently smooth initial data. 

Section \ref{derivation} sets out the fundamental field equations, which in a basic sense date back to  \ct{L}, \ct{W}.
Section \ref{periodic Cauchy problem} is devoted to local well-posedness (existence, uniqueness and continuous dependence on initial data) of  classical solutions to the  initial value problem for the pseudo-differential equation derived in Section \ref{derivation}. We establish this result by rewriting the problem on the topological group $\D$ of $C^{1}$ class diffeomorphisms as an initial value problem for an ordinary differential equation. Section \ref{global} provides   apriori estimates for  more regular classes of solutions, leading to global existence of such solutions in three and  higher dimensions. 

In Section \ref{weak}, we reconsider a class of piecewise affine solutions previously mentioned in \ct{CISY}. These solutions are less regular than those arising in our  existence results. It is found that they exist globally, independently of the underlying dimension.

\section{The {\em n\,}-Dimensional Equation}
\label{derivation}
\setcounter{equation}{0}

Consider the  {\em n}-dimensional  Euler equations for an ideal, 
 inviscid and incompressible fluid 

\ba
\dt \uu +\uu . \nabla\uu  +\nabla p& = & 0,\nonumber\\
\nabla .\uu & = & 0,\eq{2.1}
\ea
where ${\xx}=(x_1,\dots x_n) \df (x_1, {\xx'}).\,$  Denoting  $x_1$ by $x, $\,  $\uu (x, \xx', t)$ represents the 
spatial velocity field of the fluid and $p(x, \xx', t)$ 
its pressure. We impose the ansatz  
\be {\bf u}(x, \xx', t)\, =\, (u(x, t), -\dx u(x, t){\bf v}(\xx' ,t)), \eq{AA}\ee
where the $(n\!\!-\!\!1)$-dimensional vector field, $\vv$, will be chosen below.
As a consequence of \re{AA},  equation \re{2.1}\,i) may be written as
\be
\dt u + u\dx u  + \dx p = 0, \eq{2.1A}
\ee
together with
\be
\dt\dx u \,{\bf v}+\dx u\, \dt{\bf v}+u\dxx u\, {\bf v}-{(\dx u)}^2 {\bf v}\cdot\nabla '{\bf v} -
\nabla ' p ={\bf 0}, \eq{2.1AA}
\ee
where the primed operators refer to the variable $\xx'$.
Using \re{AA} and \re{2.1A}, one sees that $\nabla ' \dx p=0$. Hence,  differentiating \re{2.1AA} in $x$ eliminates the pressure term to give
\be
 \dx (\dt \dx u+u\dxx u)\, {\bf v}  +\dxx u\, \dt{\bf v}-\dx((\dx u)^{2}){\bf v}\cdot\nabla '{\bf v}={\bf 0}.
\eq{2.1AAA}
\ee
Applying the $\nabla '\cdot$\, operator to \re{2.1AAA} and using \re{2.1}\,ii) with \re{AA}, to find that $\nabla '\cdot{\bf v}=1$, shows
\be
\dx(\dt \dx u+u\dxx u)-\dx((\dx u)^{2})\, \,\nabla '{\bf v}:\nabla '{\bf v}=0,
\eq{2.1AAAA}
\ee
where $\nabla '{\bf v}:\nabla '{\bf v}={\textstyle tr}(\nabla '{\bf v})^{2}=\partial_j v_k\,\partial_k v_j$ (summing  over $j, k$ from 2 to n).
For compatibility, we must choose $\bf v$ such that $\nabla '{\bf v}:\nabla '{\bf v}$ is independent
of $\xx'$. This can be done, for instance, by choosing ${\bf v}=\frac{1}{n-1}\,\xx'$,
in which case $\nabla '{\bf v}:\nabla '{\bf v} = \frac{1}{n-1}$ and \re{2.1AA}   takes the form
\be
(\dt\dx u+u\dxx u-{\textstyle\frac{1}{n-1}}\,(\dx u)^2){\bf v} -\nabla ' p ={\bf 0}.
\eq{2.1AAAAA}
\ee

In the following we examine the periodic, initial-boundary value problem, with 
 boundary conditions for $x\in\tor \simeq \reel / \mathbb{Z} , t\geq 0,$ given by
\be
{\bf u}(0, {\bf x'}, t)= {\bf u}(1, {\bf x'}, t),\eq{2.2.2}
\ee
and
\be
p (0, {\bf x'}, t)=p (1, {\bf x'}, t). \eq{2.2.4}
\ee
Since \re{2.1AAAA} now becomes
\be
\dx(\dt\dx u+u \dxx u)-{\textstyle\frac{1}{n-1}}\,\dx((\dx u)^2))=0,\eq{2.9}
\ee
which we remark happens to be the
$x$-derivative of a Calogero-class equation, (\ct{C}),
$$\dx\dt u + u\dxx u-\Phi (\dx u) = 0,$$
 we obtain the equation

\be
\dt\dx u+u \dxx u-{\textstyle\frac{1}{n-1}\,}(\dx u)^2=f \eq{2.10}
\ee
with $f$ purely a function of time. This implies, by \re{2.1AAAAA}, 
that
\be
\nabla ' p ={\textstyle\frac{f}{n-1}\,}\,\, \xx '\eq{2.105}
\ee
 while, by \re{2.1A},
$-\dxx p={\textstyle\frac{n}{n-1}\,}(\dx u)^2+f$ 
and so
$\DD p = -{\textstyle\frac{n}{n-1}}(\dx u)^2 .$
Finally, for sufficiently smooth functions $u(x,t)$, using \re{2.2.2}  and  integrating 
 \re{2.10} we have 

\be
f = -{\textstyle\frac{n}{n-1}}\int_{\tor}(\dx u)^2 dx\eq{2.13}
\ee
while  \re{2.1A}, \re{2.2.2} and \re{2.2.4} imply that

\be
\frac{d}{dt}\int_{\tor}u\, dx = 0.\eq{2.14}
\ee




Let us introduce the  operator $\dxi$  defined by 
\[ \dxi \phi (x, t)=\int_{x_0}^{x}\phi (y, t)dy-\int_{\tor}\int_{x_0}^{x}\phi (y, t)\, dy dx .\]  
We make the observation that  
$\dx$ and $\dxi$ generally do not commute since 
$[\dx , \dxi]\phi = \int_\tor \phi\, dx, $
where $[P, Q]=PQ-QP.$ 

Consider  equation \re{2.10}, written in the form
\begin{equation}
\begin{array}{l} 
\dx (\dt u + u \dx u) = {\textstyle\frac{n}{n-1}\,} (\dx u)^{2} +f(t) 
\end{array}
\label{eq:1}
\end{equation} 
with $n>1$.
As a result of \re{2.13} and the fact that 
$\dxi\dx\phi=\phi-\int_{\tor}\phi\,dx,$
we may write equation (\ref{eq:1}) in a nonlocal form as
\begin{equation}
\dx(\dt u+u \dx u)= {\textstyle\frac{n}{n-1}\,} \dxi\dx ((\dx u)^{2})
\label{eq:2}
\end{equation}
and then, using the periodicity of $u$, we obtain
\begin{equation}
\dt u+u \dx u= {\textstyle\frac{n}{n-1}\,} \dxii\dx ((\dx u)^{2}).
\label{eq:3}
\end{equation}

Okamoto and Zhu \ct{OZ} previously established local existence  for \re{2.9}  with $u\in H^2$, using a method introduced by Kato and Lai. Their approach requires showing uniqueness separately and  then using uniqueness to prove continuous dependence on initial data. Here we instead derive a local well-posedness result in $C^{1}$ which follows from Picard iteration  after rewriting the equation as an ordinary differential equation on an infinite dimensional Banach space. This method can also be used to prove well-posedness in Sobolev spaces $\sob(\tor)$ for $s>3/2$ (see 
\ct{HM}, for example, for a similar result for the Camassa-Holm equation).

\section{Local Existence of Classical Solutions}
\label{periodic Cauchy problem}
\setcounter{equation}{0}

In  \ct{Arn}, Arnold observed that the initial value problem for the classical Euler equations of a perfect fluid
can be stated as a geometric problem of finding geodesics on the group of volume preserving
diffeomorphisms. Following this observation, Ebin and Marsden, \ct{EMa}, 
developed the functional analytic tools to establish sharp local well-posedness results for the Euler
equations. This  method has since been used for other equations with similar geometric interpretations; for example Misio\l ek,  \ct{Mis}, obtained  local well-posedness in $C^{1}(\tor)$ for the Camassa-Holm equation, which is the equation for geodesics of the $H^{1}$ metric on the Virasoro group.

In this section we develop an appropriate analytic framework for  equation  (\ref{eq:3}), using a similar approach to prove the following theorem.

\begin{theorem} 
Suppose that $n>1$.
Then there exists a unique solution 
\[ u\in C^{0}([0,T), C^{1}(\tor))\cap C^{1} ([0, T), C^0 (\tor)),
\]
to the Cauchy problem for  equation (\ref{eq:3}) with initial data $u_{0}\in C^{1}(\tor)$ for some $T>0$, and the solution depends continuously on the initial data.
\label{th:C1}
\end{theorem}

Let $ \ga$ be the flow generated by $u$, that is, $u=\gd\circ\gi$.  Then we 
obtain the equation
\begin{equation}
\gdd =  {\textstyle\frac{n}{n-1}\,} \dxii\dx \left( (\dx (\gd\circ\gi))^{2}\right) \circ\ga
\label{eq:*}
\end{equation}
from (\ref{eq:3}).
Therefore it is sufficient to prove that 
\[ F(\gd,\ga)= {\textstyle\frac{n}{n-1}\,} \left(\dxii\dx \left( (\dx (\gd\circ\gi))^{2}\right)\right)\circ\ga
\]
defines a continuously differentiable vector field in a neighborhood of the identity on the topological group $\D$ of $C^{1}$ class diffeomorphisms. Then Theorem \ref{th:C1} follows by Picard iteration  over Banach spaces. 

We remark that the smooth dependence on initial data for (\ref{eq:*}) implies only  continuous dependence on initial data for (\ref{eq:3}). The map $\ga\rightarrow \gi$ is continuous but not locally lipschitz, \ct{EMa}, and this prevents  obtaining more regularity for the initial data to solution map by this method.  The question of whether the regularity of the solution map $u_0\rightarrow u(t)$ can be improved, or not, is open. It is known, for instance, that it is not possible to improve the regularity of this map for the Camassa-Holm equation in Sobolev spaces, \ct{HM2}.
\\
In the remainder of this section, $C_\ga$ will represent  generic constant depending only on the $C^1$ norms of $\ga$ and $\gi$. 
\\

\noindent
{\bf Proof of Theorem \ref{th:C1}.} Let us denote by $P_{\ga}$ the operator given by conjugation
\[ P_{\ga}(\phi):=P(\phi\circ \gi)\circ\ga
\]
for any $\ga\in \D$ and pseudodifferential operator $P$. Using this notation we write 
\[ F(\gd,\ga)= {\textstyle\frac{n}{n-1}\,} \dxiidxg  (\dxg \gd)^{2}.
\]
Next we compute the directional derivative $\partial_{\ga}F_{(\gd,\ga)}$  and prove that it is a  bounded linear map.

Note that  $(\dxii\dx)f=\dxi \{f-\int_{0}^{1} f \ dx\}$ is a bounded operator from $C^0(\tor)$  into $C^1(\tor)$. Furthermore we abuse the notation slightly and denote by $\D$ the subgroup of orientation preserving $C^1$ diffeomorphisms which has index 2. 

Let $s\rightarrow \ga_{s}$ be a smooth curve in $\D$ such that $\ga_{0}=id$ and $\partial_{s} \ga_{s}|_{s=0}=W$ for $W \in C^{1}(\tor)$.  Then we have
\be
\partial_{\ga}F_{(\gd,\ga)}(W)=\frac{n}{n-1}\{ (\partial_{\eps}G_{\eps})|_{\eps=0} \circ \ga +W (\dx G_{\eps})_{\eps=0}\circ \ga  \}
\label{eq:10171}
\ee
where $G_{\eps}=\dxii\dx \{  (\dx(\gd\circ\gi_{\eps}))^{2} \} $.

We know that $\dx\dxi$ gives the identity, hence the second summand on the right in (\ref{eq:10171}) can be written as 
\be 
W \dx G_{\eps}|_{\eps=0}\circ\ga = \left\{ (\dxg \gd)^2 - \int_{0}^1 (\dx (\gd\circ\gi))^2 dx \right\} W 
\ee 
Moreover the computation of the first summand on the right in (\ref{eq:10171})  is reduced by
\be
 (\partial_{\eps}G_{\eps})|_{\eps=0}= \dxii\dx(\partial_{\eps}H_{\eps}|_{\eps=0})
\ee
to determine  $\partial_{\eps}H_{\eps}|_{\eps=0}$ where $H_{\eps}=(\dx(\gd\circ\gi_{\eps}))^2$. Here a straightforward computation leads to
\be
\partial_{\eps} H_{\eps}|_{\eps=0}=- \left( \dx (\gd\circ\gi) \right)^2  \dx (W\circ\gi)- \dx \{(W\circ\gi)\left(\dx(\gd\circ\gi)\right)^2 \}
\ee
Then, after an integration by parts, we obtain
\[ \begin{array}{ll}
\partial_{\eps} G_{\eps}|_{\eps=0}=&- \left(\dx(\gd\circ\gi)\right)^2 (W\circ\gi) + \int_{0}^1 \left(\dx(\gd\circ\gi)\right)^2 (W\circ\gi) dx   \\
&-  \dxii\dx \{ (\dx(\gd\circ\gi))^2  \dx(W\circ\gi)\} \\
&+ (x-\frac{1}{2}) \{ \left(\dx(\gd\circ\gi)\right)^2 (W\circ\gi)\}|_{0}^{1} 
\end{array} \]
The last term on the left hand side of the above inequality vanishes since $\gd$ and $W$ are periodic functions and $\ga$ is an orientation preserving diffeomorphism. Therefore we have
\be \begin{array}{ll}
\partial_{\ga}F_{(\gd,\ga)}(W)=\frac{n}{n-  1} \big\{&-\dxii\dx \{ (\dx(\gd\circ\gi))^2 \dx(W\circ\gi)\}\circ\ga \\ 
& +\int_{0}^1 (\dx(\gd\circ\gi))^2 W\circ\gi dx \\
&- W\int_{0}^1 (\dx(\gd\circ\gi))^2 dx \big\} 
\end{array} 
\label{eq:10191}
\ee

The linearity of the map $W\rightarrow \partial_{\ga}F_{(\gd,\ga)}(W)$ is clear. Thus we proceed to show that it is bounded. It is sufficient to estimate the $C^1$ norms of all three summands on the right in (\ref{eq:10191}). The second and third terms are both bounded by $C_\ga \|W\|_{C^0} \|\gd\|_{C^1}^2$. For the first term on the right in (\ref{eq:10191}), we have
\begin{eqnarray}
\| (\dxii\dx)_{\ga}\{ (\dxg \gd)^2 \dxg W \}\|_{C^1}  \leq & \| (\dxii\dx)_{\ga}\{ (\dxg \gd)^2 \dxg W \}\|_{C^0} \nonumber \\
& + \|  (\dxg \gd)^2 \dxg W \|_{C^0} \|\ga\|_{C^1}
\end{eqnarray}
which is bounded by $C_\ga\|\gd\|_{C^1}^2 \|W\|_{C^1}$.

In the direction of $\gd$, the G\^{a}teaux derivative of $F$ is given by
\[ \partial_{\gd}F_{(\gd,\ga)}(W)= \frac{2n}{n-1} (\dxii\dx)_{\ga} \left( \dxg \gd \dxg W \right)
\]
and this is a bounded map since
\begin{eqnarray}
\|\partial_{\gd}F_{(\gd,\ga)}(W)\|_{C^{1}} & \leq & C_{n} \| (\dxii\dx)_{\ga} \left( \dxg \gd \dxg W \right)\|_{C^{0}} \nonumber \\
& &+ \|\dxg \gd \dxg W - \int_{\tor}\dxg \gd \dxg W \ dx \|_{C^{0}} \|\ga\|_{C^{1}} \nonumber \\
&\leq & C_{n,\ga}\|\gd\|_{C^{1}} \|W\|_{C^{1}}. \nonumber
\end{eqnarray}
where $C_{n,\ga}$ depends only on $n$ and $C^1$ norms of $\ga$ and $\gi$.

In order to complete the proof of theorem \ref{th:C1} it is sufficient to show that $F$ is Fr\'{e}chet differentiable i.e, both directional derivatives $\partial_{\ga}F$ and $\partial_{\gd}F$ are continuous maps. 

{\bf \em Continuity of $(\gd, \ga) \rightarrow \partial_{\ga}F_{(\gd,\ga)}(W)$}

The following inequality reduces the proof of continuity of $\partial_{\ga}F_{(\gd,\ga)}(W)$ to estimating the two summands on the right hand side:
\be \begin{array}{ll}
\|\partial_{\ga}F_{(\gd_1,\ga_1)}(W)-\partial_{\ga}F_{(\gd_2,\ga_2)}(W)\|_{C^1}\leq & \| \partial_{\ga}F_{(\gd_1,\ga_1)}(W)-\partial_{\ga}F_{(\gd_2,\ga_1)}(W)\|_{C^1}  \\ \\
& + \| \partial_{\ga}F_{(\gd_2,\ga_1)}(W)-\partial_{\ga}F_{(\gd_2,\ga_2)}(W)\|_{C^1}
\end{array} \nonumber
\ee
where the inequality holds up to a constant depending on $n$. We rewrite the $C^1$ norm that we wish to estimate to show continuity in $\gd$ as
\begin{eqnarray}
\lefteqn{\| \partial_{\ga}F_{(\gd_1 ,\ga)}-\partial_{\ga}F_{(\gd_2,\ga)}\|_{C^{1}}} \nonumber \\
& \leq & \| (\dxii\dx)_{\ga} \left\{ ((\dxg\gd_{1})^2-(\dxg\gd_{2})^2)\dxg W \right\}\|_{C^{1}} \nonumber \\
& & + \left| \int_{0}^1 \left( \dx(\gd_1\circ\ga-\gd_2\circ\ga)\right)^2 W\circ\gi dx \right|\label{eq:bitiyo}\\
& & +  |W|\left| \int_{0}^1 \left( \dx(\gd_1\circ\ga-\gd_2\circ\ga)\right)^2  dx \right| \nonumber
\end{eqnarray} 
The last two summands in (\ref{eq:bitiyo}) are bounded by $C_\ga \| \gd_1-\gd_2\|_{C^1}^2 \|W\|_{C^0}$. For the remaining term, we have 
\be \begin{array}{ll}
\| (\dxii\dx)_{\ga} \left\{ ((\dxg\gd_{1})^2-(\dxg\gd_{2})^2)\dxg W \right\}\|_{C^{1}}  \\ \\
  \ \ \leq  \| (\dxii\dx)_{\ga} \left\{ ((\dxg\gd_{1})^2-(\dxg\gd_{2})^2)\dxg W \right\}\|_{C^{0}} \\ \\
  \ \ \ \ \ + \| \left( (\dx(\gd_1 \circ\gi))^2- (\dx(\gd_2 \circ\gi))^2 \right) \dx (W\circ\gi) \|_{C^0} \|\dx \ga\|_{C^0}
 \end{array} \nonumber
 \ee
 which is bounded by $C_\ga \| \gd_1 - \gd_2 \|_{C^1} \| \gd_1 + \gd_2 \|_{C^1} \| W\|_{C^1}$.

Our next estimate establishes continuity of $\ga\rightarrow\partial_{\ga}F_{(\gd,\ga)}(W) $. Note that it is sufficient to consider 
\begin{eqnarray}
& & \| \partial_{\ga}F_{(\gd,\ga)}(W)-\partial_{\ga}F_{(\gd,id)}(W)\|_{C^1} \nonumber \\ 
& & \ \ \leq \| (\dxii\dx)_{\ga} \left\{ (\dxg \gd)^2 \dxg W \right\} - \dxii\dx ((\dx \gd)^2 \dx W)\|_{C^1} \label{eq:100} \\ 
& & \ \ \ \ + \left| \int_0^1 \left\{ (\dx(\gd\circ\gi))^2 (W\circ\gi)-(\dx\gd)^2 W\right\}dx \right| \label{eq:101} \\ 
& & \ \ \ \ + \|W\|_{C^0} \left| \int_0^1 \left\{ (\dx(\gd\circ\gi))^2-(\dx\gd)^2 \right\}dx \right| \label{eq:102}
\end{eqnarray}
where the inequality is up to a constant depending on $n$. After adding and subtracting the appropriate terms (\ref{eq:101}) is bounded by
\be
\| W\|_{C^0} \| \dx\gd\|_{C^0} \| (\dx\gd\circ\gi) \dx\gi - \dx\gd\|_{C^0}+ \| \dx\gd\|_{C^0}^2  \|W\|_{C^1} \| \ga - id\|_{C^0}
\ee
which is bounded (up to a constant $C_{\ga}$) by $\| \gd\|_{C^1}^2  \|W\|_{C^1} \| \ga - id\|_{C^1}$. The term in (\ref{eq:102}) is estimated similarly. Hence in order to establish continuity in $\ga$ it is sufficient to bound  
\begin{eqnarray}
&\| (\dxii\dx)_{\ga} \left\{ (\dxg \gd)^2 \dxg W \right\} - \dxii\dx ((\dx \gd)^2 \dx W)\|_{C^1}  \nn \\
& \leq \| (\dxii\dx)_{\ga} \left\{ (\dxg \gd)^2 \dxg W \right\} - \dxii\dx ((\dx \gd)^2 \dx W)\|_{C^0}  \label{eq:201}\\
&+\|  \left\{ (\dx( \gd\circ\gi)^2 \dx( W\circ\gi) \right\}\circ\ga \dx\ga - (\dx \gd)^2 \dx W\|_{C^0} \label{eq:202}
\end{eqnarray}
The norm in  (\ref{eq:202}) is equal to 
\be
\| \dx W (\dx\gd)^2  \left\{ (\dx\gi)^2\circ\ga-1\right\}\|_{C_0}
\ee 
which is bounded by $C_{\ga} \|\dx\gd\|_{C^0}^2 \|\dx W\|_{C^0}\|\ga -id \|_{C^1} $. 
For (\ref{eq:201}) it is sufficient to estimate
\be
\| (\dxii\dx)_{\ga}S - \dxii\dx S\|_{C^0} + \| \dxii\dx S - \dxii\dx ((\dx\gd)^2\dx W) \|_{C^0} \label{eq:201ab}
\ee
where 
\[ S=S(\gd,\ga,W)=(\dx (\gd\circ\gi))^2\circ\ga \dx(W\circ\gi)\circ\ga .\] 
We observe that, in general, $\dxii \dx$ and $\dxi$ satisfy the identities
\ba
&&(\dxii\dx)_{\ga} f -\dxii\dx f \\
&&= \dxig f-\dxi f 
- \left( \dxi\left( \int_0^1 f\circ\gi(x)dx \right)\circ\ga - \dxi\left( \int_0^1 f(x)dx \right)\right) \nn 
\ea
and
\be
\dxig f- \dxi f = \dxi (f(\dx\ga-1))-\int_0^1 \int_{x}^{\gi(x)}f(\eta)\dx\ga(\eta , t)d\eta dx.
\ee
The second identity above is obtained changing variables in the first two integral terms in $\dxig f- \dxi f $   .  Moreover we have 
\ba
&& \dxi\left( \int_0^1 f\circ\gi(x)dx \right)\circ\ga - \dxi\left( \int_0^1 f(x)dx \right) \nn \\  
&&= \dxig \left(\int_0^1 f(\eta)(\dx\ga(\eta)-1)d\eta \right) \nn \\
&&+ \dxig \left( \int_0^1 f(x)dx \right) - \dxi \left( \int_0^1 f(x)dx\right) \nn
\ea
if both $f$ and $\dx\ga$ are periodic and $\ga$ is an orientation preserving diffeomorphism of $\tor$.  Therefore, combining these identities, we obtain the following bound for the first summand in (\ref{eq:201ab}):
\be
\|S\|_{C^0} \|\ga - id\|_{C^1} \leq C_{\ga} \|\dx\gd\|^2_{C^0} \|\dx W\|_{C^0} \|\ga -id\|_{C^1} 
\ee
 The second summand in (\ref{eq:201ab}) is equal to the $C^0$ norm of
 \be
 \dxii\dx \left( \dx W (\dx\gd)^2 \left\{(\dx\gi)^3\circ\ga-1\right\} \right)
 \ee 
 which is bounded by $C_{\ga} \| W\|_{C^1}\| \gd\|_{C^1}^2 \| \ga - id \|_{C^1}$.
 Hence the continuity of $\ga \rightarrow \partial F_{(\gd,\ga)}$ follows.

The continuity of $(\gd,\ga)\rightarrow \partial_{\gd}F_{(\gd,\ga)}$ can be shown analogously. 
Therefore $F(\gd,\ga)$ defines a continuously differentiable map in a neighborhood of $(id,u_{0})$.
This completes the proof of theorem \ref{th:C1}.
\qed

 
\section{Global Existence for $n\geq 3$}
\label{global}
\setcounter{equation}{0}

In this section we investigate the persistence of solutions of the initial value problem for (\ref{eq:3}) and show that, unlike the two-dimensional case where solutions  may blow up in finite time (\ct{CISY}, \ct{Cox}), they persist  for $n\geq 3$ in the appropriate function spaces. 

For the  theorem below let us use the following notation, 

\[ X_{n}(\tor)=\left\{ \begin{array}{ll} 
W^{2, \infty}(\tor), & n=3 \\
W^{2, \frac{n-1}{n-3}}(\tor), & n>3,
\end{array} \right. \]
and
\[ Y_{n}(\tor)=\left\{ \begin{array}{ll} 
W^{1, \infty}(\tor), & n=3 \\
W^{1, \frac{n-1}{n-3}}(\tor), & n>3,
\end{array} \right. \]

\begin{theorem} Let $n\geq 3$ and assume that $u_0(x)\in X_n(\tor)$. Then 
$u(x, t)\in C^{0}([0,T], X_n(\tor))\cap C^{1} ([0, T], Y_n(\tor))$, for arbitrary $T>0.$

\label{glob}
\end{theorem}
\noindent{\bf Proof\,}
Consider equation \re{2.9}, expressed in the form

\be
\dt\dxx u+u \dxxx u+{\textstyle\frac{n-3}{n-1}}\,\,\dx u\dxx u=0.
\eq{up2}
\ee
\noindent
Using the flow $\ga$ of $u$ ($\dot\gamma=u\circ\gamma$) in \re{up2}, we first solve for $\dx^{2} u$,
\be
\dxx u\circ\gamma(t)=u''_0\,\exp\left(-\frac{n-3}{n-1}\int_0^t\dx u\circ\gamma(s)ds\right).
\eq{reg}
\ee
By the identity $\dx\dot\gamma=\dx u\circ\gamma\,\dx\gamma$, we also have
\be
\dx\gamma(t)=\exp\left(\int_0^t\dx u\circ\gamma(s)ds\right).
\eq{alpha}
\ee
 Therefore
\be
\dxx u\circ\gamma(t)(\dx\gamma(t))^{\frac{n-3}{n-1}}=u''_0.
\eq{combine}
\ee
Note that Theorem \ref{th:C1} implies that $\gamma\in C^1$  locally in time and it follows,  given $\dx\gamma(0)~=~1$, that  there exists an interval, $t\in[0, \tau (\varepsilon))$, over which $0<\varepsilon\leq inf_{x\in\tor}\,\dx\gamma(t)\leq sup_{x\in\tor}\,\dx\gamma(t)\leq\varepsilon^{-1}.$  Equation \re{combine} then implies that, locally, $u\in X_n(\tor)$, since $\gamma$ maps $\tor$ diffeomorphically to itself.  In turn, equation (\ref{eq:2}) shows that $\dt u\in Y_n(\tor)$ over the same time interval.  (With additional assumptions on the data, further regularity can  also be bootstrapped to higher derivatives).

Assuming then that sufficient smoothness holds locally in time, we find on  multiplying \re{up2} by $|\dxx u|^{p-2}\,\dxx u$ that

\be
\dt |\dxx u|^p + u\dx |\dxx u|^p +{\textstyle p\,\frac{n-3}{n-1}}\,\dx u|\dxx u|^p=0.
\eq{up2p}
\ee
Since equations \re{AA} and \re{2.2.2} imply that both $u$ and $\dx u$ are periodic functions of $x$,  the same is true of $\dxx u$, by \re{2.10}.   One therefore obtains, on integrating \re{up2p} over $\tor$, 

\be
\frac{d}{dt}\int_\tor |\dxx u|^p dx + ({\textstyle p\,\frac{n-3}{n-1} - 1})\int_\tor \dx u|\dxx u|^p dx = 0,
\eq{int2p}
\ee 
from which it follows that the $L^\frac{n-1}{n-3}(\tor)$ norm of $\dxx u$ is uniformly conserved in time for $n>3$.
The case $n=3$ can either be considered as the limit $n\rightarrow 3$ with $p\rightarrow\infty$ in \re{int2p}, or directly using \re{up2} which  shows that $\dxx u$ is constant along characteristics and hence its $L^\infty (\tor)$ norm is uniformly conserved.

Periodicity of $u(x, t)$ in $x$ implies there exists a zero for $\dx u$, say at $x= x_0 (t),$ and so for $x, x_0\in\tor$,
 \[ \dx u (x, t)=\int_{x_0}^x \dyy u(y, t) dy.\]
For $n>3$, we  therefore  have the estimate 
\[ |\dx u (x, t)| \leq |x-x_0|^\frac{2}{n-1} ||u''_0||_{\frac{n-1}{n-3}} \leq ||u''_0||_{\frac{n-1}{n-3}}\] 
using H\"{o}lder's inequality, and so 
\[ ||\dx u||_\infty\leq||u''_0||_{\frac{n-1}{n-3}}. \]
\noindent
If $n=3$, then  
\[ ||\dxx u||_{\infty} =  ||u''_0||_{\infty}\]
which means 
\[ |\dx u (x, t)| \leq |x-x_0| ||u''_0||_{\infty} \leq ||u''_0||_{\infty} \]
for all $x~\in~\tor$, and so 
\[ ||\dx u||_\infty\leq||u''_0||_\infty \]
for all $t>0$.

Further,  since $u(x, t)-u_0(x)$ has mean zero by equation \re{2.14},  there exists $x=x_1(t)$ where  $u(x_1, t)=u_0(x_1)$ and,  for $x, x_1\in\tor$, we have
\[ u(x, t)=u_0(x)+\int_{x_1}^x (\dy u(y, t)-u'_0(y))dy.\]
It follows that
\[ |u(x, t)|\leq||u_0||_\infty + |x-x_1|(||\dx u||_\infty + ||u'_0||_\infty)\] 
for all $x~\in~\tor$, which gives the inequality
\[||u||_\infty\leq ||\dx u||_\infty + ||u_0||_{C^1}.\]

Combining the results of the previous two paragraphs shows that

\be \begin{array}{lcr}
||u||_{C^1}\leq ||u_0||_{C^2} & {\mbox for } & n=3,
\end{array} 
\eq{est2_1bottom}
\ee
and
\be\begin{array}{lcr}
||u||_{C^1}\leq ||u_0||_{C^1}+ ||u''_0||_{\frac{n-1}{n-3}} &
{\mbox for} & n>3.
\end{array}
\eq{est2_1top}
\ee

Finally, on using the properties of the operators $\dxi\dx$ and $\dxii\dx$ in equations (\ref{eq:2}) and (\ref{eq:3}) together with the above estimates, it is  seen that $||\dt u||_\infty$ and  $||\dt\dx u||_\infty$ are majorized  by a  function of $||u_0||_{C^2}$ for $n=3,$ while $||\dt u||_\infty$ and  $||\dt\dx u||_\frac{n-1}{n-3}$ are majorized by a  function of $||u_0||_{W^{2, \frac{n-1}{n-3}}}$ for $n>3.$



In both of these cases it follows  that the $C^1(\tor)$ norm of $u$ and the $C^0(\tor)$ norm of $\dt u$ remain uniformly bounded in time over any interval of local existence and, by bootstrapping the arguments of Theorem 1, the solution  can be continued, globally, in time.

\qed


As a remark, we  note here how a blow-up argument made in \ct{CISY}, which involves  a nontrivial class of separable solutions to equation \re{2.9} for $n=2$, fails to apply in the case $n>3$.  In particular, the possible appearance of a $(\tau - t)^{-1}$ factor, $\tau>0$, in the two-dimensional case  no longer exists in higher dimensions.

Given the solution form $u(x, t)=X(x)T(t),$ equation \re{2.9} reduces to
\be
\lambda X''(x)+\frac{n-3}{n-1}X'(x)X''(x) + X(x) X'''(x) = 0,\,\, x\in\tor,
\eq{X} 
\ee
where
\be
\dot{T}(t)-\lambda T(t)^2=0,\,\, t\geq 0,
\eq{T}
\ee
and $\lambda$ is a constant.
Multiplying equation \re{X} by $|X''(x)|^{\frac{5-n}{n-3}}X''(x)$ now gives
\be
\lambda|X''(x)|^\frac{n-1}{n-3}+\frac{n-3}{n-1}(X'(x)|X''(x)|^\frac{n-1}{n-3} +  X(x) (|X''(x)|^\frac{n-1}{n-3})')=0.
\eq{X''}
\ee
By using periodicity, an integration of \re{X''} over $\tor$  for $n>3$ therefore shows
\be
\lambda\int_\tor |X''(x)|^\frac{n-1}{n-3}dx = 0
\eq{Xiden}
\ee
and the result then follows. We note  that there are in general one or more points of inflection in nontrivial, periodic solutions, which prevents this argument from holding in two dimensions.


\section{Weak Solutions}
\label{weak}
\setcounter{equation}{0}
In this section, we construct a basic, piecewise differentiable class of   weak solutions $u(x, t)\in C^{0}([0,T), PC^{1}(\tor))\cap C^{1} ([0, T), PC^0 (\tor))$
to \re{2.10}, which are found to exist for all $T>0$, regardless of the underlying dimension. 

For every vector field ${\Phi}({\bf x}, t)\in C_0^\infty ([0, T)\times\tor\times\RR^{n-1};\RR^n)$ such that $\nabla\cdot{\Phi}=0,$ and for every scalar function ${\theta}({\bf x}, t)\in C_0^\infty([0, T)\times\tor\times\RR^{n-1};\RR),$ the  velocity field ${\bf u}({\bf x}, t)$ in \re{2.1} satisfies


\ba
&\int_Q \partial_t \Phi\cdot{\bf u}+ (\nabla\Phi\, {\bf u})\cdot{\bf u}\, d{\bf x}\, dt +\int_{\tor\times\RR^{n-1}} \Phi({\bf x}, 0)\cdot{\bf u}({\bf x}, 0)\, d{\bf x}= 0,
\eq{w1}
\\
&\int_{\RR^n} \nabla\theta\cdot{\bf u}\,d{\bf x}= 0,
\eq{w2}
\ea
where $Q=[0, T)\times\tor\times\RR^{n-1}.$

In terms of \re{AA}, 
equations \re{w1} and \re{w2} reduce to 

\ba
\int_Q \partial_t\phi\,u-\partial_t\Phi'\cdot{\bf v}\,\partial_x u\dbxdt
\nonumber
\\
+\int_Q \partial_x\phi\, u^2
-(\partial_x\Phi'+\nabla'\phi)\cdot{\bf v}\, u\partial_x u
+(\nabla'\Phi'\,{\bf v})\cdot{\bf v}\,(\partial_x u)^2\dbxdt
\nonumber
\\
+\int_{\tor\times\RR^{n-1}} \phi({\bf x}, 0)u(x, 0)-\Phi'({\bf x}, 0)\cdot{\bf v}({\bf x}, 0)\partial_x u(x, 0)\dbx =0,
\eq{w3}
\\
\int_{\tor\times\RR^{n-1}} \partial_x\theta\, u - \nabla'\theta\cdot{\bf v\,}\partial_x u\dbx = 0,
\eq{w4}
\ea
in which we have used the notation $\Phi=(\phi, \Phi')$ to distinguish the first component from the remaining $n-1$ components of $\Phi.$
Denoting by $[{\bf u}]={\bf u_+}-{\bf u_-}$ the jump in $\bf u$  across any smooth surface of discontinuity, $S$, and considering test functions whose support crosses $S$, equation \re{w2}  shows that 
\be
[{\bf u}]\cdot{\bf n}=0\eq{jump}
\ee
where ${\bf n}=({\mathfrak n}, {\bf n'})$ is  normal to $S$.
Then, by \re{w4}, we have 
\be
[u]{\mathfrak n} - [\partial_x u]{\bf v}\cdot{\bf n'}=0,\,\,\mbox{where\,\,\,} {\bf v}={\textstyle\frac{1}{n-1}}{\bf x'}.
\ee
In examining weak, frontlike, piecewise continuous solutions for which $[u]=0$ and $ [\partial_x u]\neq 0$ (see \ct{Daf}), it follows that these discontinuities  propagate so that ${\bf n'\cdot x'}=0.$  
A weak formulation  specific to such discontinuities may be derived by means of appropriate choice of test functions from \re{w3}, or  by
 observing that \re{2.10} may be written in conservation form as
\be
\dt((\dx u)^{1-n})+\dx(u(\dx u)^{1-n})+(n-1)(\dx u)^{-n}f=0,\,\,x\in\tor.\eq{c1}
\ee
We will admit weak solutions, $u(x, t),$ which satisfy the relation
\ba
 \int_{\mathfrak Q} \partial_t \varphi\,{(\dx u)^{1-n}}+ \dx\varphi\, u(\dx u)^{1-n}
 -(n-1)f\varphi (\dx u)^{-n}
 \, d{ x}\, dt \\
 +\int_{\tor} \varphi({x}, 0)(\dx u)^{1-n}({ x}, 0)\, d{ x}= 0,
 \ea
for all $\varphi(x, t)\in C_0^\infty (\mathfrak Q)$, where ${\mathfrak Q}=[0, T)\times\tor$.
Using  standard Rankine-Hugoniot type arguments, \ct{Daf},   discontinuities in $\dx u$ that jump across a curve $x=\psi (t)$ are seen to satisfy
\be
(-\dot\psi + u(\psi , t))[(\dx u)^{1-n}]=0 \eq{jumpdisliked}
\ee
and such discontinuities therefore propagate with the flow of \re{2.10}, {\em i.e.} $\psi (t)$ is a member of the characteristic family, $\dot\gamma =u\circ\gamma$.

\subsection{Piecewise Affine Solutions}

We begin by commenting on the general case of  periodic, $N-$phase, piecewise affine solutions. Given that both $\dx u$ and $\dt u$ may be discontinuous across the curves $x=\psi_i (t), 1\leq i\leq N-1$,  in order for $u$ to remain continuous there we must have $[u](\gamma(t),t)=0$, and so $\frac{d}{dt}[u](\gamma(t),t)=0.$ As a result, $\dt[u]+u\dx[u]=0,$ and first derivative jumps are seen to satisfy the relations

\be
[\dt u]+u[\dx u] = 0 \,\,  \mbox{and } [\dx p]=0,
\ee
 from \re{2.1A}.
Under these conditions on $u$, the expression for $f(t)$, which was obtained in \re{2.13} by integrating \re{2.10} for $u\in C^1(\tor)$, remains unchanged:
\be
f= - \frac{n}{n-1}\int_\tor (\dx u)^2 dx.
\eq{faffine}
\ee

In the special case, $N=2$, which we consider here, our form of solution  is given by the periodic extension of the function

\be
u(x, t)=\alpha +\left\{
\begin{array}{ll}x\,p, & x\in (0, \tilde{\phi} ),\\ 
\tilde{\phi} p+(x-\tilde{\phi})\,q, & x\in (\tilde{\phi} , \tilde{\psi}),\\
\tilde{\phi} p+(\tilde{\psi}-\tilde{\phi})\,q+(x-\tilde{\psi})\,p,& x\in (\tilde{\psi}, 1),
\end{array}
\right.
\eq{2phase}
\ee
where $\tilde{\phi}=\phi - [[\phi]], \tilde{\psi}=\psi-[[\psi]],$ with $[[.]]$ denoting the \lq integer part\rq\, of the argument. The functions $\tilde\phi(t)$ and $\tilde\psi(t)$ are the representatives in $[0, 1]$ of the phase curves $\phi(t), \psi(t) \in(-\infty, \infty)\,$ which start out from $ \phi(0), \psi(0)\in[0, 1]$  and  separate those regions where $\partial_x u(x, t)$ periodically takes on values of $p(t)$ or $q(t)$.

Proceeding heuristically for the moment, periodicity of $u(x, t)$  requires that
\be
{\cal N}(t)=\phi(t) p(t)\ +(\psi(t)-\phi(t))q(t)+(1-\psi(t))p(t)=0.
\eq{phipsi}
\ee
Given the spatial periodicity in pressure (see \re{2.2.4}), we recall that integration  of equation \re{2.1A} over one period showed the integral $\int_0^1 u(x, t)dx$ to be independent of time. This allows \re{2phase} to be used to give an expression for $\alpha (t).$ The result may be written, for instance in terms of $p,  \tilde\phi$ and $\tilde\psi$, as
\be
\alpha + \frac{p}{2}(\tilde\phi + \tilde\psi -1) = c
\eq{ap}
\ee
where we have set  $c=\int_0^1 u(x, 0)\, dx=\alpha(0)+\frac{p(0)-q(0)}{2}(\psi(0)-\phi(0))(\phi(0)+\psi(0)-1)$. We will assume that $0<\phi(0)<\psi(0)<1$. Choosing, for convenience, $c=0$, the \lq average characteristic\rq\, must propagate with speed zero and we will see, consequently, that both $\phi (t)$ and $\psi(t)$ remain in $[0, 1]$ for all time. The distinctions between $\tilde\phi$ and $\phi$, $\tilde\psi$ and $\psi$, will therefore not be made further here.

Combining \re{2.10}, \re{faffine} and \re{2phase}  now gives
\be
\dot{p}=\frac{1}{n-1}p^2+f,\,\,\dot{q}=\frac{1}{n-1}q^2+f
\eq{pandq}
\ee
where
\be
f=-\frac{n}{n-1}(\phi p^2+(\psi-\phi)q^2+(1-\psi)p^2).
\eq{f}
\ee
Also, by \re{jumpdisliked},
\be
\dot\phi=\alpha+\phi p
\eq{620}
\ee
and
\be
\dot\psi
=\alpha+\phi p+(\psi-\phi)q.
\eq{621}
\ee

We next verify that \re{phipsi} follows from the system of equations \re{pandq} - \re{621}.
Differentiation and some simplification gives
$$
\dot{\cal N}=-(\psi - \phi)(p-q)q + (\psi - \phi)q^2+(\phi + (1-\psi))(\frac{p^2}{n-1}+f)+(\psi-\phi)(\frac{q^2}{n-1}+f)
$$
$$
=-(\psi - \phi)(p-q)q-((\phi+(1-\psi))p^2+(\psi-\phi)q^2)
$$
$$
=-p(q(\psi - \phi)+p(\phi+(1-\psi))
$$
$$
=-p{\cal N},
$$
and so ${\cal N}(t)={\cal N}(0)e^{-\int_0^t p(s)ds}$. In particular, taking periodic data, ${\cal N}(0)=0,$ means that ${\cal N}(t)=0,\, t>0.$
We may therefore use \re{phipsi} to write \re{621} as
\be
(1-\psi)\dot{}=-\alpha+(1-\psi)p.
\eq{622}
\ee
and again employ \re{phipsi}  to express the following phase fractions as functions of $p$ and $q$,
\be
\psi-\phi=\frac{p}{p-q}\,, \,\phi + (1-\psi) = \frac{-q}{p-q}\,.
\eq{psiphi}
\ee
Using these relations in \re{f} leads to 
\be
f=\frac{n}{n-1}pq
\ee
 from which \re{pandq} reduces to 
 the autonomous system
\ba
\dot{p}=\frac{1}{n-1}(p^2+n\,pq),
\eq{22i}\\
\dot{q}=\frac{1}{n-1}(q^2+n\,pq). 
\eq{22ii}
\ea

Subtracting \re{620} from \re{621} implies
\be
\psi (t)-\phi (t)=(\psi (0)-\phi (0))\exp(\int_0^t q(s)ds)
\eq{623}
\ee 
which means that the center phase fraction does not collapse as long as $\int_0^t q(s)ds$ remains bounded away from $-\infty$. 
Similarly, adding \re{620} and \re{622} gives
\be
\phi (t)+(1-\psi (t))=(\phi (0)+(1-\psi (0))\exp(\int_0^t p(s)ds)
\eq{624}
\ee
and the outer phase fraction exists as long as $\int_0^t p(s)ds>-\infty$.
Comparing equations \re{623} and \re{624} shows also that periodicity imposes the following requirement on $p(t)$ and $q(t)$ in terms of their initial phase fractions,
\be
(\phi (0) + (1-\psi (0))\exp(\int_0^t p(s)ds) + (\psi (0)-\phi (0))\exp(\int_0^t q(s)ds)=1.
\eq{phases}
\ee

Using \re{pandq} to compute $p-q$ next gives
\be
p(t)-q(t)=(p(0)-q(0))\exp(\frac{1}{n-1}\int_0^t p(s)+q(s)\, ds),
\eq{p-q}
\ee
 which shows $p(t)-q(t)$ cannot change sign. By equation \re{psiphi}, $p(t)$ and $q(t)$ consequently keep their signs as long as neither phase fraction collapses.  This can alternatively be seen by considering a sketch of $u$  and observing that $p$ and $q$ have opposite signs and can only vanish simultaneously.
Thus, without loss of generality, we  subsequently assume $p(t)>0>q(t)$, for at least some $t\geq 0$.

In the following lemma, we solve the system \re{22i}, \re{22ii} implicitly in order to show that two-phase solutions of the type \re{2phase} exist for all time.
 
\begin{lemma}{Let $(p(t), q(t))$ satisfy \re{22i}, \re{22ii} with initial data $p(0)>0>q(0)\,\, ($respectively, $p(0)<0<q(0))$.}
Then $(p(t), q(t))$ exists for all $t\in(-\infty, \infty)$
and satisfies $p(t)>0>q(t)\,\, ($respectively, 
$p(t)<0<q(t))$. Further,\, $(p(t), q(t))\rightarrow (0, 0)$ as $t\rightarrow\pm\infty$ and $||u(., t)||_{C^1}+||\dt u(., t)||_{C^0}\rightarrow 0$ as $t\rightarrow\pm\infty.$
\end{lemma}
\noindent{\bf Proof} Writing equations \re{22i} and \re{22ii} using polar variables, \linebreak $p(t)=r(t)\cos \theta (t),\,\, q(t)=r(t)\sin \theta (t), r(t)\geq 0,$  leads to the following,
\be
\dot{r}(t)=\frac{r^2}{n-1}(\cos^3(\theta(t))+n(\ctt+\stt)\ctt\stt+\sin^3(\theta(t))),
\eq{rdot}
\ee
and
\be
\dot{\theta}(t)=r(\ctt-\stt)\ctt\stt.
\eq{thetadot}
\ee
In the case  $p(0)>0>q(0)$,  $\theta(0)\in(-\pi/2 , 0)$, so by \re{thetadot} $\dot{\theta} (0)<0$
and, as long as $\theta (t)\in (-\pi/2, 0),$ $\dot{\theta}(t)<0$. 

We will show that $\theta(t)\rightarrow -\pi/2$ as $t\rightarrow\infty$
($\theta(t)\rightarrow 0$ as $t\rightarrow-\infty$) by using \re{rdot} and 
\re{thetadot}. Integrating the resulting expression for $\frac{d\ln r}{d\theta}$  gives

\be
r(\theta)=c|\cos\theta\sin\theta|^{\frac{1}{n-1}}|\cos\theta-\sin\theta|^{-\frac{n+1}{n-1}}
\eq{rtheta}
\ee
where $c> 0$  denotes a generic constant. Inserting this expression for $r(\theta)$ in \re{thetadot} results in
\ba
\dot{\theta}(t)=c|\cos\theta\sin\theta|^\frac{1}{n-1}|\cos\theta-\sin\theta|^{-\frac{n+1}{n-1}}(\cos\theta-\sin\theta)\cos\theta\sin\theta\nonumber\\
= -c|\cos\theta\sin\theta|^\frac{n}{n-1}|\cos\theta-\sin\theta|^{-\frac{2}{n-1}}\, ,\hspace{22 ex}
\eq{theta}
\ea
for $\theta (t)\in (-\pi/2, 0)$, and so 

\be
\int_{\theta(0)}^{\theta(t)}\frac{|\cos\theta-\sin\theta|^\frac{2}{n-1}}{|\cos\theta\sin\theta|^\frac{n}{n-1}}d\theta=-ct,\,\, -\frac{\pi}{2}<\theta(0)<0,\,\, c>0.
\eq{int}
\ee
Since $\frac{n}{n-1}>1,$ the integral expression diverges, both as $\theta(t)\rightarrow-\pi/2$ \,($t\rightarrow+\infty)$, and as $\theta(t)\rightarrow0$ ($t\rightarrow -\infty).$ Thus $\theta (t)$ and, from \re{rtheta}, $r(t),$ are bounded, continuous functions of time, with $\theta (t)\in (-\pi/2, 0)$.  By \re{rtheta} then, $r(t)\rightarrow 0$ as $t\rightarrow\pm\infty$ 
and, noting that  from \re{phipsi}, \re{ap},
$$\alpha (t)=-\frac{p(t)-q(t)}{2}(\psi(t)-\phi(t))(\phi(0) + \psi(0) -1)
\, ,$$
it follows from  \re{indivphi} and \re{indivpsi} that 
$$p(t), q(t), \alpha(t) \mbox{ and } u(x, t)\rightarrow 0 \mbox{ as } t\rightarrow\pm\infty.$$
The remaing conclusions are  easily obtained.\qed
\linebreak

Finally, we show that for $c=0$ the phases remain in the interval $[0, 1]$ for all time.

\begin{theorem} Suppose $c=0$ and $0<\phi(0)<\psi(0)<1$. Then the phases $\phi, \psi$ stay in $[0, 1]$. In particular $$\phi(t)\in[\phi(0), \frac{1}{2}(\phi(0)+\psi(0)) \mbox{\, and\,\,\,} \psi(t)\in(\frac{1}{2}(\phi(0)+\psi(0)), \psi(0)]$$ for all $t>0.$  Further, the time asymptotic behavior satisfies $$\lim_{t\rightarrow\infty}\phi(t)=\frac{1}{2}(\phi(0)+\psi(0))=\lim_{t\rightarrow\infty}\psi(t)$$.
\end{theorem}
\noindent{\bf Proof}

With \re{ap}  giving $\alpha$ (for $c=0$), substituting into \re{620} and \re{622}  shows that, by \re{624}, 
\be
\dot\phi(t)=(1-\psi(t))\dot{}=\frac{p}{2}(\phi+(1-\psi))=\frac{1}{2}(\phi(0)+(1-\psi(0)))p(t)\exp(\int_0^t p(s)ds)
\ee
and the individual phases therefore satisfy
\be
\phi(t)=\phi(0)+\frac{1}{2}(\phi(0)+(1-\psi(0)))(\exp(\int_0^t p(s)ds)-1)
\eq{indivphi}
\ee
and
\be
\psi(t)=\psi(0) -\frac{1}{2}(\phi(0)+(1-\psi(0)))(\exp(\int_0^t p(s)ds)-1).
\eq{indivpsi}
\ee

Now we examine \re{phases}. Assuming $q<0<p$, the first term is monotonically increasing in time and must converge, for $0<\phi(0)<\psi(0)<1$, to a positive limit. The second term is positive but monotonically decreasing, so may converge, as $t\rightarrow\infty$ either to a positive limit, or to zero.
Formally, setting $0<\int_0^\infty p(t)dt=L<\infty,$ and $ -\infty\leq \int_0^\infty q(t)dt=M<0,$ equations \re{phases} and \re{p-q} imply

\be
(\phi(0)+(1-\psi(0)))e^L+(\psi(0)-\phi(0))e^M=1
\eq{infty1}
\ee
and
\be
lim_{t\rightarrow\infty}(p(t)-q(t))=(p(0)-q(0))e^\frac{L+M}{n+1}.
\eq{infty2}
\ee
If $M$is finite and $p(0)\neq q(0)$, then \re{infty2}, together with the necessity for
$p(t)$ to approach zero as $t\rightarrow\infty$,  means that $\lim_{t\rightarrow\infty} q(t)\neq0.$ However this implies that $M=-\infty$, a contradiction. On the other hand, if $M=-\infty$ then, since $L<\infty$, \re{infty2} requires that $\lim_{t\rightarrow\infty}q(t)=0$, which is permitted. We conclude that, for smooth solutions $p(t), q(t)$ to exist with $p(0)>0>q(0)$, we require that $L=\int_0^\infty p(t)dt<\infty$ and $M=\int_0^\infty q(t)dt=-\infty.$ Thus, by \re{infty1},
\be
e^L=\frac{1}{\phi(0)+1-\psi(0)}.
\eq{L}
\ee
Writing \re{indivphi} in the form
\be
\phi(t)=\frac{1}{2}(\phi(0)-(1-\psi(0))+\frac{1}{2}(\phi(0)+(1-\psi(0)))\exp\int_0^t p(s)ds
\eq{indivphi1}
\ee
it follows that
\be
\lim_{t\rightarrow\infty}\phi(t)=\frac{1}{2}(\phi(0)+\psi(0))
\eq{phiinfty}
\ee
and similarly,
\be
\lim_{t\rightarrow\infty}\psi(t)=\frac{1}{2}(\phi(0)+\psi(0)).
\eq{psiinfty}
\ee
By the monoticity in time of $\int_0^t p(s)ds$, therefore
$\phi(t)\in[\phi(0), \frac{1}{2}(\phi(0)+\psi(0))$ and $\psi(t)\in(\frac{1}{2}(\phi(0)+\psi(0)), \psi(0)]$ for all $t>0.$


\end{document}